\documentclass[12pt]{article} 
\usepackage[utf8]{inputenc} 

\usepackage{geometry}
\usepackage{graphicx}
\usepackage{tikz}
\usetikzlibrary{trees,shapes,arrows}

\usepackage{paralist}
\usepackage{subfig}
\usepackage{sidecap}
\usepackage{array}

\usepackage{amsthm}
\usepackage{amsmath}
\usepackage{amsfonts}
\usepackage{mathtools}

\usepackage{float}
\usepackage{hyperref}
\hypersetup{
	colorlinks=false,
	linkbordercolor=gray,
	citebordercolor=gray,
	urlbordercolor=gray} 
\usepackage{cleveref} 

\usepackage{sectsty}
\allsectionsfont{\sffamily\mdseries\scshape\centering\nohang}

\usepackage{abstract}

\usepackage{caption}
\usepackage{sidecap}

\allowdisplaybreaks

\usepackage[switch,modulo]{lineno}
\usepackage{amssymb}
\usepackage[affil-it]{authblk}

\usepackage{todonotes}
\usepackage{lineno}

\newtheoremstyle{mplain}
  {2\topsep}  
  {\topsep}  
  {\itshape} 
  {0pt}      
  {\scshape}
  {:}        
  {5pt plus 1pt minus 1pt}
  {}

\renewenvironment{proof}[1][\proofname]{{\scshape #1. }}{\qed}

\theoremstyle{mplain}

\newtheorem{thm}{Theorem}[]
\crefname{thm}{Theorem}{Theorems} 
\Crefname{thm}{Theorem}{Theorems}

\newtheorem{lemma}[thm]{Lemma}
\crefname{lemma}{Lemma}{Lemmas} 
\Crefname{lemma}{Lemma}{Lemmas}

\newtheorem{conj}[thm]{Conjecture}
\crefname{conj}{Conjecture}{Conjectures} 
\Crefname{conj}{Conjecture}{Conjectures}

\newtheorem{remark}[thm]{Remark}
\crefname{remark}{Remark}{Remarks} 
\Crefname{remark}{Remark}{Remarks}

\newtheorem*{claim}{Claim}

\newtheorem{prop}[thm]{Proposition}
\crefname{prop}{Proposition}{Propositions} 
\Crefname{prop}{Proposition}{Propositions}

\crefname{cor}{Corollary}{Corollaries} 
\Crefname{cor}{Corollary}{Corollaries}

\newtheorem{definition}[thm]{Definition}
\crefname{definition}{Definition}{Definitions} 
\Crefname{definition}{Definition}{Definition}

\newcommand{\E}{\mathbb{E}}
\newcommand{\R}{\mathbb{R}}
\newcommand{\N}{\mathbb{N}}
\newcommand{\C}{\mathcal{C}}
\renewcommand{\Pr}{\mathbb{P}}
\newcommand{\eps}{\varepsilon}

\newcommand{\lint }{\int_{\Lambda}}
\newcommand{\fcn}{\psi}
\DeclareMathOperator{\poi}{Poi}

\DeclareMathOperator{\Wei}{Weibull}

\providecommand{\keywords}[1]{\textbf{\textit{Keywords---}} {\small #1}} 
\begin{document}
%\linenumbers
\title{\huge The Minimum Perfect Matching \\ in Pseudo-dimension $0<q<1$}
\author{Joel Larsson
\thanks{Electronic address: \texttt{joel.larsson@warwick.ac.uk}}}
\affil{Mathematics Institute,\\University of Warwick, United Kingdom}
%\date{January 19, 2019}

\maketitle
\begin{abstract}
It is known that for $K_{n,n}$ equipped with i.i.d.\ $\exp(1)$ edge costs, the minimum total cost of a perfect matching converges to $\pi^2/6$ in probability. Similar convergence has been established for all edge cost distributions of \emph{pseudo-dimension} $q \geq 1$, for instance $\Wei(1,q)$ or $\chi_q$ (chi distribution with $q$ degrees of freedom). In this paper we extend those results to all real positive $q$, confirming the Mézard-Parisi conjecture in the last remaining applicable case.
\end{abstract}
\keywords{matching, mean field, replica symmetry, random graph, pseudo-dimension}

\label{intro}
\newpage
\section{Introduction}
\label{background}

There has been substantial interest over the past few decades in the minimum matching problem: Given a graph $G$, and a positive cost associated to each edge of $G$, we want to find a perfect matching of minimal total cost $M(G)$. 
Of special interest is minimum matching on the complete bipartite graph $K_{n,n}$ on ${n+n}$ vertices with random edge costs given by independent $\exp(1)$-variables, sometimes referred to as the \emph{random assignment} problem. 
For this graph model, the lower bound $\liminf_n M(K_{n,n})\geq 1$ of minimum matching is trivial (the cheapest edge from any given vertex has expected cost $n^{-1}$, and a perfect matching uses $n$ edges).
The upper bound $\limsup_n M(K_{n,n}) \leq 3$ was established by Walkup \cite{walkup2}, by finding a perfect matching using only fairly cheap edges. This was later improved to $2$ by Karp \cite{karp-upperbound}.

Mézard \& Parisi \cite{repopt} conjectured that $M(K_{n,n})$ converges in probability to $\zeta(2)=\pi^2/6$, based on heuristic replica symmetry calculations. Aldous \cite{asymptotics} proved that the limit exists, and later confirmed the conjecture \cite{zeta2}. Both of these papers used what is sometimes called the `objective method', and worked with matchings on an infinite limit object.
Parisi \cite{parisiconj} further conjectured the more precise result that $\E[M(K_{n,n})] = \sum_{k=1}^n k^{-2}$. This was later established independently by Nair, Prabhakar \& Shaw \cite{nairetal} and Linusson \& Wästlund \cite{linussonwastlund}, both using inductive proofs. The proof was later simplified by Wästlund \cite{wastlundinductivezeta2}.
Salez and Shah \cite{salezshah} gave yet another proof of the Mézard-Parisi conjecture, using the objective method to analyze the behavior of belief propagation on the limit object. 

A more comprehensive overview of the existing literature and related problems can be found in a survey paper by Krokhmal and Pardalos \cite{krokpar}.

A natural question is whether these results extend to other edge cost distributions. It turns out that only the scaling behaviour of the probability distribution near $0$ matter.
A random graph where edge costs are i.i.d.\ copies of the random variable $\ell$ is said to be of \emph{pseudo-dimension} $q$ if $\lim_{x\to 0} \Pr(\ell \leq x) \cdot x^{-q}$ exists and lies in $(0,\infty)$. 
The distributions $\operatorname{Exp}(1)$ and $\operatorname{Unif}(0,1)$ are both of pseudo-dimension $1$, and the chi distribution with $q$ degrees of freedom as well as the Weibull distribution with shape parameter $q$  (i.e.\ the $(1/q)$:th power of an exponential variable) are of pseudo-dimension $q$.
Mézard \& Parisi \cite{repopt} considered these distributions for real positive $q$, but most focus since then has been on the special case $q=1$.

The motivation for the term pseudo-dimension is this: For $q \in \N$, a geometric graph model is given by embedding the vertices as $n$ points chosen uniformly and independently at random in a hypercube $[0,1]^q$, and setting the edge costs to be the corresponding Euclidean distances. The mean field approximation (i.e.\ the graph model where edge costs follow the same distribution, but are i.i.d.) of a geometric graph model of dimension $q$ is a graph model of pseudo-dimension $q$.

For any graph $G$ and probability measure $\nu$ on $\R_+$, let $G[\nu]$ denote $G$ equipped with i.i.d edge costs with distribution given by $\nu$. If $\nu$ is of pseudo-dimension $q$, then the cost of the minimum matching on $K_{n,n}[\nu]$ can be shown to be of order $n^{1-1/q}$, by a minor modification of \cite{walkup2}. This suggests studying the quantity $n^{-1+1/q}M(K_{n,n}[\nu])$. Does it converge in probability to a constant?
This question was answered in the affirmative for $q\geq 1$ by Wästlund \cite{repsym} (both for $K_{n,n}$ and $K_n$), but it remained open for $0<q<1$. 
Our main result is the following theorem, confirming the Mézard-Parisi conjecture for all $q>0$.
\begin{thm}
\label{mainthm}
For every $q>0$, there exists a $\beta=\beta(q)$ such that
for any probability measure $\nu$ for which $c:= \lim_{x\to 0}\nu(\{\ell\leq x\})\cdot  x^{-q}$ exists and $c\in(0,\infty)$,
\[
\frac{M(K_{n,n}[\nu])}{c^{1/q}n^{1-1/q}} \to \beta(q)
\quad\textrm{and}\quad
\frac{M(K_{2n}[\nu])}{c^{1/q}(2n)^{1-1/q}} \to \frac{1}{2}\beta(q)
\]
in probability as $n\to \infty$.
\end{thm} 
We believe the theorem should hold in somewhat greater generality, i.e.\ for graphs $G$ other than $K_{n,n}$ and $K_n$. As we will discuss in \cref{relaxedmatchings-explorationgame}, the important properties of $G$ are its local structure (tree-like when expensive edges are removed) and its expansion properties. We therefore conjecture the following.  
\begin{conj}
\label{graphonconjecture}
Let $G_n$ be a sequence of graphs on $2n$ vertices which admits a graphon limit $\mathcal{G}$. Let $W$ be the function on $[0,1]^2$ given by $\mathcal{G}$.
If there is a $\delta>0$ such that $\int_0^1 W(x,y)dx=\delta$ for all $y \in [0,1]$, then $n^{-1+1/q}M(G_n[\nu])$ converges in probability to a constant (depending only on $\delta$, $q=q(\nu)$,  and $c=c(\nu)$).
\end{conj}

\label{prelim} 
\section{Notation and definitions}

We will assume $0<q<1$ and a large parameter $\lambda>1$ is fixed, and often suppress dependence on them in our notation. 
We will restrict our attention to $0<q<1$, since \cref{mainthm} is already known to be true for $q\geq 1$. 
Although our proof strategy works for $q\geq 1$ too, some parts of our lemmas are trivial in that case, and assuming $0<q<1$ streamlines the proofs. 
Unless otherwise stated, all functions considered will be real-valued functions on $\Lambda:=[-\lambda/2,\lambda/2] $. For $f$ and $g$ functions on $\Lambda$, we will use ${f\leq g}$ to mean that ${f(z)\leq g(z)}$ for all $z\in\Lambda$. For $x\in \R$, we let \(x_+ := \max(x,0)\).
We write $a_n\nearrow a$ if $a_n\to a$ and $a_n$ is a non-decreasing sequence, and $a_n\searrow a$ if $-a_n\nearrow -a$. 
For a weighted rooted graph $G$, the $(k,\lambda)$-truncation $G(k,\lambda)$ is the subgraph of $G$ obtained as the $k$-neighbourhood of the root after all edges of weight more than $\lambda$ have been removed. Equivalently, $G(k,\lambda)$ is the union of all paths from the root of length at most $k$ that only use edges of weight at most $\lambda$.

We say that the rooted tree $T$ is the \emph{$m$-rooted $\lambda$-local limit} of the graph sequence $G_n$ iff for every $k\in  \N$ and every choice of $m$ roots in $G_n$, $G_n(k,\lambda)$ converges to $m$ disjoint independent copies of $T(k,\lambda)$ in the total variation metric. In the special case $m=1$, we simply say that $T$ is the \emph{$\lambda$-local limit} of $G_n$. 

Furthermore, we will let $|G|$ denote the number of \emph{edges} of a graph $G$
and we will consider the edges of a rooted tree to be directed away from the root $\phi$. By \emph{path} we will mean a directed path away from the root. If $u$ is the parent of $v$, we write $u \to v$. Let $|u|$ denote the distance to $u$ from the root.
For any edge $uv$ we let $\ell(u,v)$ denote its cost. 

\section{Proof strategy}
As a first step, we switch to working with a rescaled model: Multiply all edge weights in $K_n[\nu]$ and $K_{n,n}[\nu]$ by $n^{1/q}$. Since minimum perfect matching is a linear programming problem, the only effect this has on the optimum is to multiply it by the same amount. Rescaling the problem in this way allows us to study $\lambda$-local limits of these graphs more easily. One can think of this as changing the units of costs in such a way that the expected number of edges of cost at most $1$ from a given vertex is $1$.

Let $\tilde K_n$ and $\tilde K_{n,n}$ be these rescaled models (suppressing dependence on $\nu$). \Cref{mainthm} is then equivalent to existence of $\lim_n n^{-1}M(\tilde K_n)$ and $\lim_n  n^{-1}M(\tilde K_{n,n})$.

\subsection{Exploration Game} 

The game \emph{Exploration} was introduced in \cite{repsym}.
This zero-sum, perfect information game is played in the following way: On an edge-weighted rooted graph $G$,
Alice and Bob takes turns picking the next edge of a self-avoiding walk starting from the root. When it is a player's turn (Alice's, say), and the current vertex is $u$, she can take one of two actions:
\begin{enumerate}[(i)]
\item Pick any neighbour $v$ of $u$ that has not already been visited, and pay Bob the cost $\ell(u,v)$ of the edge $uv$. Bob then continues the game from $v$.

\item Quit the game, and pay Bob a penalty of $\lambda/2$, for some fixed parameter $\lambda>0$.
\end{enumerate}
The payoff for Alice, once the game has finished, is the total amount Bob has payed to her minus the total amount she has payed to Bob. Each player's aim is to maximize their payoff.   
If the weighted graph $G$ is finite, every game position has a well-defined game value $f=f(G,u)$.\footnote{By the game value of $u \in V(G)$ we will mean the value of the exploration game on $G$, starting from $u$, to the second player. In other words, the net amount that Alice will pay to Bob, assuming optimal play by both.}
If Alice starts by moving from $u$ to $v$, the remainder of the game is equivalent to a game played on $G-u$ started on $v$, but with the roles of Alice and Bob reversed.
By considering all possible options a player has from the vertex $u$, it is easy to see that
\begin{equation}
f(G,u) = 
\min\Big(\lambda/2,\min\limits_{v} \big(\ell(u,v) - f(G-u,v)\big),
\Big),
\label{gamevaluerecursiongeneral}
\end{equation}
where the second minimum is taken over all neighbours $v$ of $u$. If the graph is a finite tree $T$ and we start the game at the root, no move can go from a vertex to its parent, so we may as well forbid such moves. But then $f(T,v)=f(T-u,v)$ if $u$ is the parent of $v$, so we let $f(v):=f(T,v)$. Thus \cref{gamevaluerecursiongeneral} can be slightly simplified to
\begin{equation}
f(u) = 
\min\Big(\lambda/2,\min\limits_{v} \big(\ell(u,v) - f(v)\big)
\Big),
\label{gamevaluerecursion}
\end{equation}
where the second minimum is taken over all \emph{children} $v$ of $u$.
If the tree is infinite, however, it is no longer clear that the function $f$ is well-defined. 
Instead, we consider all functions $f$ which satisfy \cref{gamevaluerecursion} for all $u\in V(T)$, and call these `game valuations'. 

It was proven in \cite{repsym} that for any $q>0$, the limits of $n^{-1}M(\tilde K_n)$ and $n^{-1}M(\tilde K_{n,n})$ exist \emph{if} there for all large $\lambda$ exists a unique game valuation on (almost all realizations of) a certain edge-weighted Galton-Watson tree $T^q_\lambda$, which we will define shortly. They proceeded to prove the valuation was indeed unique for $q\geq 1$ (Proposition 2.8 in \cite{repsym}), but that proof did not extend to $0<q<1$. Therefore, in order to prove \cref{mainthm} it suffices to show the following: 
\begin{prop}
\label{prop:unique valuation}
For any $\lambda>0$ and $q\in (0,1)$ there is almost surely a unique game valuation on $T^q_\lambda$, i.e.\ a function $f:V(T^q_\lambda)\to \Lambda$ satisfying \cref{gamevaluerecursion} for every $u\in V(T^q_\lambda)$.
\end{prop}

The tree $T^q_\lambda$ is constructed in the following way.
For every vertex $u$, run an inhomogeneous Poisson point process on $[0,\lambda]$ with intensity $qt^{q-1}$ at time $t$. (Call this intensity measure $m_\ell$.) If $\ell_1,\ell_2,\ldots \ell_k$ are the times of the events in this process, we let $v_1,v_2,\ldots v_k$ be the children of $u$, and give the edge $uv_i$ weight $\ell_i$. By Proposition 2.2 of \cite{repsym}, $T^q_\lambda$ is the $m$-rooted $\lambda$-local limit of $\tilde K_n$ and of $\tilde K_{n,n}$, for any $m\in \N$.

An equivalent definition is that $T^q_\lambda$ is the Galton-Watson tree with offspring distribution $\poi(\lambda^q)$, and i.i.d.\ edge weights given by the probability measure $m_\ell/\lambda^q$. 
This tree is also related to the `Poisson-weighted infinite tree' (PWIT). Let $T$ be the graph obtained by raising all edge weights in the PWIT to the $(1/q)$:th power. Then $T^q_\lambda$ has the same distribution as $T(\infty,\lambda)$, i.e.\ the connected component of the root after all edges of weight more than $\lambda$ are removed.  
 
The recursion (\ref{gamevaluerecursion}) has a useful monotonicity property: if $f$ and $g$ are game valuations such that $f(v)\leq g(v)$ for all children $v$ of $u$, then $f(u)\geq g(u)$. This gives rise to the following lattice ordering:
\[
f\preceq g \iff
\begin{cases}
f(u)\leq g(u), \textrm{ $\forall u:|u|$ is odd }
\\
f(u)\geq g(u), \textrm{ $\forall u:|u|$ is even}
\label{latticeorder}
\end{cases}
\]
Let $f^k(u)$ be the unique game valuation on $V(T^q_\lambda(k,\lambda))$ (i.e.\ the tree $T^q_\lambda$ truncatad after $k$ generations) satisfying \cref{gamevaluerecursion} for all $u$ with $|u|\leq k$. Note that $f^k(u)=\lambda/2$ for all $u$ with $|u|=k$, since these $u$ have no offspring in $T^q_\lambda(k,\lambda)$.  
The following claim follows from the monotonicity of the recursion (\ref{gamevaluerecursion}).
\begin{claim} The lattice order $\preceq$ is bounded, with unique maximum and minimum given by the pointwise limits
\begin{align}
f_A:= \lim_{k\to \infty} f^{2k+1}
\quad \textrm{ and }\quad
f_B:= &\lim_{k\to \infty} f^{2k}. \label{claim:sandwich}
\end{align}
\end{claim}
\noindent To prove \cref{prop:unique valuation} it suffices to show that $f_A=f_B$ (\cite[p.1072]{repsym}).  
\subsection{Proof strategy for  \cref{prop:unique valuation} for $q\geq 1$}
\label{wastlundsproof}
We give here a short description of Wästlund's proof of \cref{prop:unique valuation} for $q\geq 1$, in order to explain how our proof for $q>0$ is similar to, yet differs significantly from, it.
Roughly speaking, there are two main components of the proof: (\ref{prob-to-quit-bound}), show that a game where both Alice and Bob plays according to $f_A$ must finish after finitely many moves, and  (\ref{few-reasonable-moves}) show this game is not `too different' from a game where Alice plays according to $f_A$ and Bob according to $f_B$.

Let $u_0:=\phi,u_1,u_2,\ldots$ be the (finite or infinite) game path when both Alice and Bob play according to $f_A$. Let $Z_i := f_A(u_i)$ for $i\geq 0$, and note that these random variables are not independent.  If $Z_i = \lambda/2$, then the $f_A$-optimal move from $u_i$ is to quit and pay the penalty, and the game path is finite iff this happens for some $i$. It's not too hard to show that ${\Pr(Z_i = \lambda/2)>0}$ uniformly in $i$, but since the $Z_i$'s are not independent this is not sufficient. However, conditional on $Z_i$, $(Z_{0},\ldots ,Z_{i-1})$ and $(Z_{i+1},Z_{i+2},\ldots)$ are independent. So in order to prove (\ref{prob-to-quit-bound}), it therefore suffices to show that there is an $\eps>0$ such that $\Pr(Z_{i+1} = \lambda/2|Z_i=z)>\eps$ for all $z\in \Lambda$. 
A move is said to be $\delta$-\emph{reasonable} iff it is within $\delta$ of being $f_A$-optimal, i.e.\ the move $u\to v$ is $\delta$-reasonable if $ \ell(u,v)-f_A(v)\leq  f_A(u) + \delta$.
It turns out that for any $\delta>0$, if Bob plays according to $f_B$ his moves will be $\delta$-reasonable eventually.
Using \cref{gamevaluerecursion} and the choice of $f_A$, Wästlund proved the following more precise statements of (\ref{prob-to-quit-bound}) and (\ref{few-reasonable-moves}):
Pick $u\in V(T)$ in a way that doesn't depend on the subtree rooted in $u$. Then
\begin{enumerate}[(i)]
\item there is an $\eps>0$ such that if $u\to v$ is the $f_A$-optimal from $u$, then\label{prob-to-quit-bound}
\[\inf_{z\in \Lambda}\Pr\big(f_A(v)=\lambda/2\big| f_A(u)=z\big)>\eps\]
\item for any $\eps>0$, there exists a $\delta>0$ such that
\[\sup_{z\in \Lambda}\E[\# (\delta\textrm{-reasonable moves from }u)|f_A(u)=z] <\eps/2.\] \label{few-reasonable-moves}

\end{enumerate}
Assume that Alice and Bob plays according to $f_A$ and $f_B$ respectively, and let $v$ be a vertex along the game path such that Bob's moves are $\delta$-reasonable after $v$. Together (\ref{prob-to-quit-bound}) and (\ref{few-reasonable-moves}) imply that the tree of all $\delta$-reasonable game paths starting from $v$ has branching number at most $1-\eps/2$, whence it is a.s.\ finite. This tree is guaranteed to contain the game path from $v$ onwards, and thus the game must end after finitely many moves. This, in turn, implies that $f_A=f_B$. 

\subsection{Revised proof strategy for $0<q<1$}
The main trouble that arises when trying to generalize the argument above to all $q>0$ is that the proof of the statement (\ref{prob-to-quit-bound}) above fails for $q<1$. Indeed, the statement is false for $q\leq \frac{1}{2}$, see \cref{q-one-half}. 

Our aim is still to show that the probability of a player quitting at any given time is uniformly bounded away from $1$ (when both players play according to $f_A$). Whether or not a player will quit the game at $u_i$ is determined by the random variable $Z_i:=f_A(u_i)$, but as mentioned earlier these random variables are not independent. This problem was sidestepped in Wästlund's proof by the conditioning in (\ref{prob-to-quit-bound}) above, but since that fails for small $q$ we will instead need to understand the dependency between $Z_i$ and $Z_{i+1}$. We do this by constructing a pair of linear operators (one for each parity of $i$) that map functions of the form ${z\mapsto \E[\bullet |Z_{i+1}=z]}$ to ${z\mapsto  \E[\bullet |Z_{i}=z]}$. The statement that will correspond to (\ref{prob-to-quit-bound}) will be that the composition of these two operators is a contraction (\cref{normbound}). 
Having changed one major component of the proof, the second one (\ref{few-reasonable-moves}) is no longer compatible. The linear operators we construct can only provide information about the conditional expectation of random variables, so we must change our aim from proving that the game path is almost surely finite, to proving that it has finite expected size. However, the expected size of the tree of reasonable moves (as defined in \cref{wastlundsproof}) does not appear to be finite.

We solve this by using a refined concept of reasonable moves, where we take into account not only single deviations by Bob from $f_A$, but instead consider the \emph{sum} of these deviations along a game path.
This leads to a significantly smaller tree of reasonable moves (guaranteed to contain the game path) whose expected size we can bound recursively.  
\subsection{The connection between Exploration and Matching}
\label{relaxedmatchings-explorationgame}
Here we will briefly discuss the connection between the seemingly disparate topics of the exploration game and the minimum perfect matching problem. This is not strictly necessary in order to understand our proof of \cref{prop:unique valuation} (which is the novel result of this paper), but it gives some insight into how \cref{prop:unique valuation} implies \cref{mainthm}. For a full proof of this implication, we refer the interested reader to section 3 of \cite{repsym}. Crucially, at no point in the proof of \cite[Thm 3.2]{repsym} is the assumption $q\geq 1$ used, the proof only depends on there being a unique game valuation $f=f_A=f_B$ on $T_\lambda^q$.

We begin by defining the \emph{$\lambda$-relaxed} (or \emph{$\lambda$-diluted}) matching problem. 
A \emph{partial matching} in a graph $G$ is a subgraph $H\subseteq G$ where no two edges share a vertex.
For a partial matching $H$ on an edge-weighted graph $G$, we say that the \emph{$\lambda$-relaxed} cost $c_\lambda(H)$ of it is the sum of the costs of all edges it contains, plus $\lambda/2$ for each unmatched vertex. In other words, $c_\lambda(H):= \sum_{uv\in E(H)}\ell(uv)+\sum_{u\notin V(H)}\lambda/2$.
We let $M_\lambda(G):=\min_H c_\lambda(H)$, where the minimum is taken over all partial matchings. Note that for any graph $G$, $M_\lambda(G)$ is an increasing function of $\lambda$, and ${M_\lambda(G)\leq M(G)}$. 
There are two results from \cite{repsym} that connects $\lambda$-relaxed matchings to perfect matchings and to the exploration game, respectively.

First, \cite[Thm 3.2]{repsym} shows that the existence of a unique game valuation (\cref{prop:unique valuation}) implies that $\beta_\lambda:= \lim_n n^{-1}M_{\lambda}(\tilde K_n)$ exists for any $\lambda$.
A rough outline of the argument is as follows.

It can be shown that for finite graphs, an arbitrary edge $uv$ participates in the optimal $\lambda$-relaxed matching iff the move $u\to v$ is optimal in the exploration game on $\tilde K_n$ starting from $u$.
It follows from \cref{gamevaluerecursiongeneral} that this move is optimal iff $\ell(u,v)\leq f(\tilde K_n,u)+f(\tilde K_n-u,v)$. The $\lambda$-local limit of both $\tilde K_n$ and $\tilde K_n-u$ is $T_\lambda^q$, and by the sandwiching argument in \cref{claim:sandwich} together with \cref{prop:unique valuation}, $f(\tilde K_n,u)$ and $f(\tilde K_n-u,v)$ can be well approximated by $f_A(\phi)$.
Letting $Z,Z'$ be two i.i.d.\ copies of $f_A(\phi)$ and $X$ be the random variable on $[0,\lambda]$ with probability density function $t\mapsto qt^{q-1}/\lambda^{q}$, the expected contribution of $uv$ is within a factor $1+o(1)$ of
\begin{align*}
\E [ \chi_{\{\ell(u,v)\leq Z+Z'\}}\cdot \ell(u,v)] &= \Pr(\ell(u,v)\leq \lambda)\cdot  \E [ \chi_{\{\ell(u,v)\leq Z+Z'\}}\cdot \ell(u,v)|\ell(u,v)\leq \lambda]
\\
&= \frac{1+o(1)}{n}\cdot \underbrace{ \lambda^q\cdot \E [ \chi_{\{X\leq Z+Z'\}}\cdot X]}_{=:\beta_\lambda},
\end{align*}
where $\chi_E$ denotes the indicator random variable for the event $E$.
Summing over all edges gives $\E M_\lambda(\tilde K_n)$, and a second moment argument shows concentration of $M_\lambda(\tilde K_n)$ around its mean.
Hence $M(\tilde K_n)\geq M_\lambda(\tilde K_n)= (1+o(1))\binom{n}{2}\frac{1}{n}\beta_\lambda$. 

This argument works $\tilde K_{n,n}$ just as well as for $\tilde K_n$, and furthermore for graph sequences $G_n$ with similar local structure. More precisely, it works for $G_n$ if the $4$-rooted $\lambda$-local limit of $G_n$ is $T_\lambda^q$. 
Second, \cite[Prop 3.4]{repsym} uses a variation on Posa's extension-rotation method to show that if the graph sequence $G_n$ satisfies an expansion property (which both $K_n$ and $K_{n,n}$ do satisfy), then a partial matching can be extended to a perfect matching at a small extra cost. Hence $M(G_n) \leq (1+o(1))M_\lambda(G_n)$.
\Cref{mainthm} follows with $\beta:=\lim_{\lambda\to \infty} \beta_\lambda$.

To motivate \cref{graphonconjecture}, note that in order to prove \cref{mainthm} for some graph sequence $G_n$, one only needs to prove two things: That the $4$-rooted $\lambda$-local limit of $G_n$ is $T_\lambda^q$, and that $G_n$ is a sufficiently good expander. It seems likely that the graphon in \cref{graphonconjecture} satisfies these properties. Some other graphs that might be good candidates are quasi-random graphs, $k$-partite complete graphs, and Erdős-Rényi graphs.   
\section{Proof of main theorem}

\subsection{The tree $T^q_\lambda$ conditional on the game valuation $f_A$} 

In order to be able to construct the linear operators mentioned above, we will change slightly how we generate the random tree $T^q_\lambda$ and the game valuation $f_A$. Instead of first generating $T^q_\lambda$ and then calculating $f_A$ `back from infinity', we will generate the tree and vertex labels $f_A$ concurrently. This will require the following lemma from \cite{repsym}.

Let $F_A(z):=\Pr(f_A(\phi)\geq z)$.
In a slight abuse of notation, we will also use $F_A$ to refer to the probability measure on $\Lambda$ of the random variable $f_A(\phi)$. Similarly, $F_B$ will refer to both the function $z\mapsto\Pr(f_B(\phi)\geq z)$ and the corresponding measure. Recall that $m_\ell$ is the measure on $[0,\lambda]$ for which $dm_\ell(t) = qt^{q-1}dt$.
Let the $\ell f$-square be the set $\{(\ell,f): 0\leq\ell\leq \lambda, |f|\leq \lambda/2\}$.
\begin{lemma}[Lemma 2.6 of \cite{repsym}\footnote{The lemma in  \cite{repsym} only states that $(\ell_i,f_i)$ constitutes a Poisson point process, not what the intensity measure is. However, it is implicit in the proof of the lemma that $\mu_B$ is the correct measure when $|u|$ is even, and the other case is analoguous.}]
\label{2d-Poisson}
Let $u\in V(T^q_\lambda)$, let $v_1,v_2,\ldots v_k$ be its children, let ${\ell_i:=\ell(u,v_i)}$, and let ${f_i := f_A(v_i)}$.
Then the points $(\ell_i,f_i)$ constitute a two-dimensional inhomogeneous Poisson point process on the $\ell f$-square, with intensity given by measure $\mu_A:=m_\ell \times F_A$ if $|u|$ is odd and $\mu_B:=m_\ell \times F_B$ if $|u|$ is even. 
\end{lemma}
An immediate consequence of the lemma is that the $f_A$-optimal move from a vertex is a.s.\ unique:  $\ell-f$ has continuous distribution, because its probability density function is given by the convolution of the function $t\mapsto qt^{q-1}$ and the measure $dF_A(-z)$.

To generate the tree $T_\lambda^q$ concurrent with $f_A$, start by picking $z$ according to the probability measure $F_A$, and assigning the root $\phi$ the game value $f_A(\phi)=z$. Then, we generate the next generation of the tree by the Poisson point process of \cref{2d-Poisson}, conditioned on ${\min(\lambda/2,\min_i(\ell_i-f_i))=z}$.

If $z=\lambda/2$, this is equivalent to conditioning on there being no point in the region $\{\ell-f< \lambda/2\}$. Since the distribution of points in two disjoint regions are independent, the points in $\{\ell-f\geq  \lambda/2\}$ are generated by an inhomogenous Poisson point process according to the measure $\mu_B$ restricted to the region ${\{\ell-f\geq  \lambda/2\}}$.

If $z<\lambda/2$, this is equivalent to conditioning on there being no point in the region $\{\ell-f< z\}$, and one special point on the line $\{\ell-f=z\}$. The points $\{\ell-f\geq  z\}$ can be generated by restricting the intensity measure to $\{\ell-f\geq z\}$. The line $\{\ell-f= z\}$ has zero $\mu_B$-measure, so to pick a random point from it we condition on there being at least one point in the Poisson point process on the strip $\{z\leq \ell-f \leq  z+\eps\}$, and then let $\eps\to 0$. Since $\mu_B = m_\ell \times F_B$, and $m_\ell$ is absolutely continuous with respect to the Lebesgue measure, this is well-defined. In order to express the probability measure obtained in the limit explicitly we must first understand the measures $F_A$ and $F_B$ in more detail. 

The following is proved in \cite[p.1077]{repsym} (as well as occuring in similar forms in e.g.\ \cite{zeta2,repopt,salezshah}), but we include the proof here because it helps in understanding some of our argument later on, in \cref{normbound,lemma:conditionalshift}.
\begin{lemma}
\label{semi-fixpoint}
Let $V$ be the non-linear operator on functions on $\Lambda$ defined by
\[
V(G)(z):=\exp\Big(-\lint q(z+t)_+^{q-1}G(t)dt \Big).
\]
Then $F_A=V(F_B)$ and $F_B=V(F_A)$.
\end{lemma} 
\noindent \begin{proof}
Recall that $F_A(z):=\Pr(f_A(\phi)\geq z)$. Since (by definition) $f_A(\phi)=  {\min(\lambda/2,\min_i(\ell_i-f_i) )}$, the event $\{f_A(\phi)\geq z\}$ happens iff there is no $(\ell_i,f_i)$ with $\ell_i-f_i<z$. By \cref{2d-Poisson}, the $(\ell_i,f_i)$ constitutes a Poisson point process, and the probability that no $(\ell_i,f_i)$ falls in the set $D_z:= \{(\ell,f):\ell-f<z\}$ is $\exp(-\mu_B(D_z))$. To calculate $\mu_B(D_z)$, first fix $\ell$ and let $t$ be such that $z+t=\ell$. Then $\ell-f<z$ iff $f>t$.
Integrating over all $t$ gives
\[\mu_B(D_z)= \lint q(z+t)^{q-1}_+\Pr(f>t)dt.\]
However, $\Pr(f>t)=\Pr(f\geq t)$ for all but countably many $t$, and $\Pr(f\geq t)=F_B(t)$ by definition.
Hence $F_A(z) =\exp(-\mu_B(D_z)) = V(F_B)(z)$. The other case is analogous.
\end{proof}

The operator $V$ is the composition of an integral operator (with a continuous kernel) and a smooth function (applied pointwise), and we can use this to establish smoothness properties of $F_A$ and $F_B$, as well as  bound their derivatives. 
\begin{lemma}
\label{locallipschitz}
Each of the two measures given by $F_A$ and $F_B$ on $\Lambda$ is the sum of a point mass at $\lambda/2$ and a measure that is absolutely continuous with respect to the Lebesgue measure on $(-\lambda/2,\lambda/2)$. The functions $F_A$ and $F_B$ are continuously differentiable on $(-\lambda/2,\lambda/2)$, with derivative $F'_A$ given by
\begin{equation}
F'_A(z) = -F_A(z)\cdot \Bigg(F_B(\lambda/2)\, q(z+\lambda/2)^{q-1}-\lint  q(z+t)_+^{q-1}F'_B(t)dt \Bigg).
\label{def:FA'}
\end{equation}
Furthermore, $ \lint q(z+t)_+^{q-1}F'_A(t)dt$ is a continuous function of $z$, and for some constant $\alpha>1$ and all ${|z|< \lambda/2}$, we have the bounds
\begin{align}
-F'_A(z) &\leq \alpha (\lambda/2-|z|)^{q-1} ,
\label{bound:F'}
\\
-\lint q(z+t)_+^{q-1}F'_A(t)dt &\leq \alpha \max\big((z+\lambda/2)^{2q-1},|z-\lambda/2|^{q-1} \big).
\label{bound:diag}
\end{align}
\Cref{bound:diag} also holds for $\lambda/2<z\leq 3\lambda/2$, and analogous results hold for $F'_B$.
\end{lemma}
\noindent The proof of this lemma is largely a lengthy calculus exercise, and we postpone it to the end of the paper. 
We will often parametrize the diagonal line $\{(\ell,f):\ell-f=z\}$ as $\{(z+t,t):t\in \Lambda\}$. The measure $\mu_B$ has  density
\begin{equation}
\rho_B^z(t) := 
q(z+t)_+^{q-1}\cdot (-F_B'(t))
\label{def:measure}
\end{equation}
 along such a diagonal for $t<\lambda/2$, and a point mass $q(z+\lambda/2)^{q-1}F_B(\lambda/2)$ at the end point $t=\lambda/2$. (And analogously for $\mu_A, \rho_A^z$.) 

\begin{lemma}
\label{diagonalmeasure}
Let $(\ell,f)$ be a point in the inhomogeneous Poisson point process on the $\ell f$-square with intensity measure $\mu_A$, conditioned to lie on the line $\ell-f=z$ (for some $z<\lambda/2$). Then the probability distribution of $f$ is given by
\[
\Pr(f < x)
=
 \frac{\int_{-\lambda/2}^x \rho_A^z(t)dt}{J^z_A}
,\quad\quad\quad
\Pr(f = \lambda/2) = \frac{q(z+\lambda/2)^{q-1}F_A(\lambda/2)}{J^z_A},
\]
where $J^z_A := q(z+\lambda/2)^{q-1}F_A(\lambda/2) +\lint \rho_A^z(t)dt$.
\end{lemma}
\noindent \begin{proof}
Let $\eta_{\eps}$ be defined as $\eps^{-1}$ times the measure $\mu_A$, restricted to the region $E_\eps:= {\{(\ell,f):z\leq \ell-f \leq z+\eps, f\geq -z\}}$, and let $\eta$ be the measure on $E_0$ which is given by $\rho_A^z(t)dt$ at the point $(z,z+t)$, and a point mass of $F_A(\lambda/2)$ at $(z,z+\lambda/2)$.

We will show that $\eta_{\eps}\to \eta$, as $\eps\to 0$, and that normalizing $\eta$ gives the probability measure in the statement of the lemma. For $z<\lambda/2$,
\[
\iint_{ f <x} d\eta_\eps 
=
\int_{-\lambda/2}^x\int_{f+z}^{f+z+\eps} \frac{1}{\eps} \rho_A^{\ell-f}(f)d\ell df
=
 - \int_{-\lambda/2}^x\frac{1}{\eps}\int_{f+z}^{f+z+\eps}  q\ell^{q-1}d\ell f_A'(f) df
\]
Note that $\ell\mapsto q\ell^{q-1}$ is a decreasing function, whence ${\frac{1}{\eps} \int_{f+z}^{f+z+\eps}q\ell^{q-1}d\ell }\nearrow {q(z+f)^{q-1}}$.
By the monotone convergence theorem, $\iint_{ f <x} d\eta_\eps \to \int_{-\lambda/2}^x  \rho_A^z(t)dt$. Similarly, the $\eta_\eps$-measure of the line segment $\{f=\lambda/2, z\leq \ell-f \leq z+\eps\}$ approaches $q(z+\lambda/2)^{q-1}F_A(\lambda/2)$ as $\eps\to 0$. So $\eta_\eps\to \eta$ as $\eps\to 0$, and 
${J_A^z := \|\eta\|\in(0,\infty)}$. (For a measure $m$, $\|m\|$ denotes the $m$-measure of the whole space on which $m$ is defined.) It follows that $\eta_\eps/\|\eta_\eps\|\to \eta/\|\eta\|$, and $\eta_\eps/\|\eta_\eps\|$ is the probability measure for a random point picked according to $\mu_A$ in $E_\eps$.  
\end{proof}

\noindent We will also use ineq. (\ref{bound:diag}) of \cref{locallipschitz} in another (weaker) form, as a bound on the normalizing factor $J^z_A$.
\begin{align}
J_A^z &\leq  \alpha \max\big((z+\lambda/2)^{2q-1},|z-\lambda/2|^{q-1} \big)+ F_A(\lambda/2)q(z+\lambda/2)^{q-1} \nonumber
\\
&<\alpha\lambda^q \max\big((z+\lambda/2)^{q-1},|z-\lambda/2|^{q-1} \big).
\label{bound:Jz}
\end{align}
\subsection{$(u,t)$-reasonable moves}
We will now introduce our new definition of reasonable moves, and show that the game path is reasonable according to this definition.

For $v\in T^q_\lambda- \{\phi\}$, let $\delta(v)$ be how  far from $f_A$-optimal it is to move to $v$ from its parent $u$. More precisely, ${\delta(v) := \ell(u,v) -f_A(u)-f_A(v)}$. Note that $\delta(v)\geq 0$, since it follows from (\ref{gamevaluerecursion}) that $f_A(u) \leq \ell(u,v) - f_A(v)$ for any $v$.

\begin{definition}
We say that a (finite or infinite) path $P=u,u_1,u_2,\ldots$ away from the root is \emph{$(u,t)$-reasonable} if $\sum_{i=1}^{|P|} \delta(u_i) \leq t$ and $\delta(u_i) = 0$ whenever $|u_i|$ is odd.
\end{definition}
\noindent In other words, a path is $(u,t)$-reasonable if Alice's moves are $f_A$-optimal and Bob's deviations from $f_A$ sum to at most $t$.
\begin{lemma}
The game path (when Alice plays according to $f_A$ and Bob according to $f_B$) is $(\phi,2\lambda)$-reasonable.
\label{reasonable-game-path}
\end{lemma}
\label{gamelength}
\label{game}
\noindent\begin{proof}[Proof of \cref{reasonable-game-path}]
Let $P$ be the game path. Pick any length $2$ sub-path $(u \to v \to w) \subseteq P$ , such that $u$ is at even distance from $\phi$.

Since $u$ is at even distance from the root, it will be Alice's turn to move from $u$. She will choose the $f_A$-optimal move, i.e.\ she will move to a child $v$ of $u$ such that $f_A(u) = \ell(u,v) -f_A(v)$. In other words, $\delta(v) = 0$. This move may or may not be $f_B$-optimal, but $f_B(u) \leq \ell(u,v) -f_B(v)$ regardless. Thus\footnote{A similar argument is used in \cite[p.1076]{repsym} to show that the difference $f_A(u_{2k})-f_B(u_{2k})$ is monotone in $k$.} 
\[
f_A(u)-f_B(u) \geq  \big(\ell(u,v) -f_A(v)\big) - \big(\ell (u,v)-f_B(v)\big) = f_B(v) - f_A(v).
\]
Then it will be Bob's turn to move from $v$. He will choose the $f_B$-optimal move, i.e.\ he will move to a child $w$ of $v$ such that $f_B(v) = \ell(v,w) -f_B(v)$. This may or may not be the $f_A$-optimal move, but by the definition of $\delta$ we have that $f_A(v) = \ell(v,w) -f_A(w)-\delta(w)$. Thus
\[
f_B(v)-f_A(v) = [\ell(v,w) -f_B(w)] - [\ell(v,w) -f_A(w)-\delta(w)]= f_A(w) - f_B(w) + \delta(w),
\]
and together with the move $u\to v$ this gives that
\[
f_A(u)-f_B(u) \geq f_A(w) - f_B(w) + \delta(w) = f_A(w) - f_B(w) +\delta(v)+ \delta(w)
\]
Let $\phi=u_0 \to u_1  \to u_2 \to  \ldots $ be the game path $P$. Pick $n\in\N$ such that $2n\leq |P|$ ($P$ might be infinite, in which case we just pick any $n\in\N$). If we repeat the argument above with $(u,v,w) := (u_{2i-2},u_{2i-1},u_{2i})$, for all $1\leq i \leq n$, we get that 
\begin{align*}
f_A(\phi)-f_B(\phi)
&\geq f_A(u_2)-f_B(u_2) + \delta(u_1)+\delta(u_2)
\\
&\geq f_A(u_4)-f_B(u_4) + \delta(u_1)+\delta(u_2)+\delta(u_3)+\delta(u_4)
\\
&\;\;\vdots
\\
&\geq f_A(u_{2n})-f_B(u_{2n}) + \sum_{i=1}^{2n}\delta(u_i)
\end{align*}
Since $|f_A|,|f_B| \leq  \lambda /2$, this implies that $\sum_{i=1}^{2n}\delta(u_i) \leq 2\lambda$.
Recall that $\delta(u_i)=0$ for odd $i$, whence $\sum_{i=1}^{2n+1}\delta(u_i) \leq 2\lambda$ as well.
So $\sum_{i=1}^{k}\delta(u_i) \leq 2\lambda$ for any $k\leq |P|$, and taking the supremum
over such $k$ it follows that ${\sum_{i=1}^{|P|}\delta(u_i) \leq 2\lambda}$.
\end{proof}

Let $\Delta_t(u)$ be the union of all $(u,t)$-reasonable paths. The crucial property of $\Delta_t(u)$ is that the event $\{w\in \Delta_t(u)\}$ is determined by the first $|w|$ generations $T_\lambda^q$, i.e.\ it is independent from the descendants of $w$. The same cannot be said of the event $\{w\in P\}$, where $P$ is the game path. The latter event depends on $f_B(w)$ (at least if $|w|$ is even), and $f_B(w)$ might not be independent from the descendants of $w$, not even conditional on $f_A(w)$. But our aim is to bound $\E|P|$, and since $P\subseteq \Delta_{2\lambda}(\phi)$ it suffices to bound $\E|\Delta_{2\lambda}(\phi)|$. 
  
We will work with $k$-level truncations $\Delta^k_t(u):=\Delta_t(u)(k,\lambda)$, and recursively bound the expected value of $|\Delta^k_t(u)|$.
Conditioned on $f_A(u)$, the distribution of $\Delta^k_t(u)$ is the same for every $u$ at even distance from the root, so we let
\begin{equation}
R^k_t(z) := \E\big[ |\Delta^k_t(\phi)| \big| f_A(\phi) = z \big].
\label{def:Rk}
\end{equation}

\begin{prop}
\label{finitegame}
There exists a family of continuous functions $(\fcn_t)_{t\in[0,2\lambda]}$ on $\Lambda$ s.t. ${R^{2k}_t(z)<\fcn_t(z)}$ for all $z \in \Lambda $, $t\in [0,2\lambda]$, and $k\in \N$, and satisfying $\sup_{z,t}\fcn_t(z)<\infty$. In particular, $\E |\Delta_{2\lambda}|<\sup_z\fcn_{2\lambda}(z)$ is finite.
\end{prop}
\subsection{Linear operators}
\noindent To prove \cref{finitegame} we will need the following lemmas concerning certain linear operators.
These operators relate functions of the form $z\mapsto \E[\bullet|f_A(u)=z]$ to $z\mapsto \E[\bullet|f_A(v)=z]$, whenever $u\to v$ is an $f_A$-optimal move.

Recall that $J^z_A := q(z+\lambda/2)^{q-1}F_A(\lambda/2) -\lint q(z+t)_+^{q-1}\cdot F_A'(t)dt$ is the measure of the diagonal line $\{\ell-f=z\}$ in the $\ell f$-square, according to the measure from \cref{diagonalmeasure}.
\vspace{-1em}
\begin{lemma}
\label{lemma:conditionalshift}
\label{def:TA}
Let the positive linear operator $L_A$ on $\mathcal{C}(\Lambda )$ and the function $I_A:\Lambda\to[0,1]$ be defined by
\begin{align}
L_A h(z) &:= \frac{\lint h(t)\rho_A^z(t)dt}{J^z_A} \label{def:T}
\\
I_A(z) &:= \frac{\lint \rho_A^z(t)dt}{J^z_A} \label{def:I}
\end{align}
on $(-\lambda/2,\lambda/2)$, and by their continuous extensions at $\pm \lambda/2$. Let also $L_B$ and $I_B$ be defined similarly.
Let $u,v$ be such that $\phi \to u \to v$ are $f_A$-optimal moves.
Then the following holds:
\begin{align}
(L_B \circ L_A) R^k_t(z) &= \E\big[ |\Delta^k_t(v)| \big| f_A(\phi) = z \big]. \label{TTonR}
\end{align}
Furthermore, $I_A$ satisfies these properties:
(i) $I_A$ is continuous, (ii) $I_A(z)<1$ for ${z \in [-\lambda/2,\lambda/2)}$, and (iii) $I_A(\pm \lambda/2)$ are well defined by continuous extension.
Analogous statements hold for $I_B$.
\end{lemma} 
\begin{lemma}
\label{normbound}
$\| L_B\circ L_A\|<1$, where $\|\cdot\|$ is the operator norm given by the $\infty$-norm on $\C(\Lambda)$.
\end{lemma} 
 
\noindent \begin{proof}[Proof of \cref{lemma:conditionalshift}]
\label{proof:conditionalshift}
Assume that the moves $\phi\to u$ and $u\to v$ are $f_A$-optimal. 
Let $Z_z := (f_A(v)|f_A(u)=z)$, and consider ${\E\big[ |\Delta^k_t(v)| \big| f_A(u) = z \big]}$. Since (by definition) ${R^k_t(z)=\E\big[ |\Delta^k_t(v)| \big| f_A(v) = z \big]}$, we can write 
\[
\E\big[ |\Delta^k_t(v)| \big| f_A(u) = z \big]
=
\E R^k_t(Z_z)
\]
For any fixed $z$, the mapping  $R^k_t\mapsto \E R^k_t(Z_z)$ is a linear functional, so the function $z\mapsto \E R^k_t(Z_z)$ depends linearly on the function $R^k_t$. But what is the linear operator that takes $R^k_t$ to $z\mapsto \E R^k_t(Z_z)$?
The distribution of $f_A(v)$ conditional on $f_A(u)=z$ is given by \cref{diagonalmeasure}. Integrating over $\Lambda$ gives that
\[
 \E R^k_t(Z_z) = 
\frac{q(z+\lambda/2)^{q-1}F_A(\lambda/2)R_t^k(\lambda/2)+\lint R_t^k(s)\rho_A^z(s)ds}{J_A^z} 
\]
But $R^k(\lambda/2)=0$, since Alice's optimal move from a vertex with game value $\lambda/2$ will be to quit immediately. Thus $ \E R^k_t(Z_z)= {\lint R^k_t(f)\rho_A^z(f)df}/{J_A^z} $, which equals $L_A R^k_t(z)$ by \cref{def:T}.
Note also that $\rho^z_A(t),\rho^z_B(t),J^z_A$ and $J^z_B$ are positive for all $z,t$, so the operators $L_A,L_B$ are positive.
Applying the same method one more time gives the desired result for the first part of the lemma. For the second part, we verify that (i)-(iii) hold.
\begin{enumerate}[(i)]
\item The non-negative term $\lint \rho_A^z(t)dt$ is continuous in $z$ by \cref{locallipschitz}, and so is the positive term ${q(z+\lambda/2)^{q-1}F_A(\lambda/2)}$. Hence both numerator and denominator of \cref{def:I} are continuous, and the denominator is non-zero, so $I_A$ is continuous.

\item Both  ${q(z+\lambda/2)^{q-1}F_A(\lambda/2)}$ and $\lint \rho_A^z(t)dt$ are positive and finite for $|z|<\lambda/2$, so $I_A(z)<1$ for such $z$.

\item Using \cref{bound:diag}, we see that for $z$ near $-\lambda/2$,
\begin{align*}
I_A(z) &= \frac{O\big((z+\lambda/2)^{2q-1}\big)}{(z+\lambda/2)^{q-1}+O\big((z+\lambda/2)^{2q-1}\big)}
=O\big((z+\lambda)^q\big),
\shortintertext{so that $\lim_{z \to -\lambda/2} I_A(z) =0$. Near $\lambda/2$,}
I_A(z) &= \frac{ \lint \rho_A^z(t)dt}{q\lambda^{q-1}F_A(\lambda/2)+o(1)+ \lint \rho_A^z(t)dt} 
\\
&= 1 -\left({1+\frac{\lint \rho_A^z(t)dt}{q\lambda^{q-1}F_A(\lambda/2)+o(1)}}\right)^{-1},
\end{align*}
so $\lim_{z\to\lambda/2} I_A(z)$ will exist if $\lim_{z\to\lambda/2}\lint \rho_A^z(t)dt$ exists (even if the latter limit is infinite).

Note that $\lint \rho_A^z(t)dt = \lint \rho_A^z(t-z+\lambda/2)dt$, as the support of $\rho_A^z$ is ${[-z,\lambda/2]\subseteq \Lambda}$, and translating by  $-z+\lambda/2$ gives a function with support ${[-\lambda/2,z]\subseteq \Lambda}$.
By \cref{bound:F'}, we have that
\[ \begin{split}
\rho_A^{z}(t-z+\lambda/2) &= q(t+\lambda/2)^{q-1}\cdot F'_A(t-z+\lambda/2)
\\
&\leq \begin{cases} \alpha q(t+\lambda/2)^{2q-2}, &t\leq 0 \\ K, &t>0, \end{cases}
\end{split} \]
for some constant $K$ and all $z$ sufficiently close to $\lambda/2$. Thus we have an upper bound on $\rho_A^{z}(t-z+\lambda/2)$ which is independent of $z$. For $q > 1/2$, this upper bound is integrable and by dominated convergence it follows that
\[
\lim_{z\to\lambda/2}\lint \rho_A^{z}(t-z+\lambda/2)dt = \lint \lim_{z\to\lambda/2} \rho_A^{z}(t-z+\lambda/2)dt = \lint \rho_A^{\lambda/2}(t)dt<\infty.
\]
Hence $\lim_{z\to\lambda/2}\lint \rho_A^z(t)dt$ exists (and is finite) for $q>1/2$. For $q\leq 1/2$, we use 
\cref{def:FA'} of \cref{locallipschitz} to replace $F'_A$:
\[
\lint \rho_A^z(t)dt \geq \lint  q(t+\lambda/2)^{q-1} F_A(t) \cdot \big(F_B(\lambda/2) \cdot q(t+\lambda/2)^{q-1}\big) dt.
\]
This integral goes to $\infty$ as $z\to\lambda/2$, since the singularity $(t+\lambda/2)^{2q-2}$ is not integrable. We conclude that $\lim_{z\to\lambda/2} \lint \rho_A^z(t)dt$ exists for all $q$, hence $I_A(\lambda/2)$ is well defined.
\end{enumerate}
\vspace{-2.07em}
\end{proof}

\begin{remark}
\label{q-one-half}
It follows from the proof of the previous lemma that $I_A(\lambda/2)=1$ for $q\leq \frac{1}{2}$, while $I_A(\lambda/2)<1$ for $q> \frac{1}{2}$. This implies that the statement (\ref{prob-to-quit-bound}) in \cref{wastlundsproof} is true only if $q>\frac{1}{2}$.
\end{remark} 
\noindent \begin{proof}[Proof of \cref{normbound}]
$L_A$ is a substochastic operator\footnote{A positive linear operator $T$ given by $T(h)(z) := \lint h(t)\kappa^z(t)dt$ is said to be \emph{stochastic} if $\lint \kappa^z(t) dt=1$ for every $z$ and \emph{substochastic} if $\lint \kappa^z(t) dt\leq 1$ for every $z$.}, and to be able to fully leverage this property we will factorize it into a stochastic operator that has almost all the structure of $L_A$ and a substochastic operator that is also a diagonal map.
Start by defining the kernel $\kappa_A^z(t)$, as $\rho_A^z$ normalized for ${(z,t) \in (-\lambda/2,\lambda/2)^2}$:
\begin{equation}
\kappa_A^z(t) := \frac{\rho_A^z(t)}{\lint \rho_A^z(s)ds}. \label{def:dk}
\end{equation}
Using this kernel, we write $L_A(h)(z)$ as $\lint I_A(z)h(t)\kappa_A^z(t)dt$. The factor $I_A(z)$ does not depend on $t$, so it can be factored out of the integral. We can therefore write $L_A$ as the composition of the operators $S_A$ and $D_A$, defined by
\begin{align}
S_A(h) (z) &:= \lint  h(t) \kappa_A^z(t)  dt
\\
D_A(h)(z) &:= I_A(z) \cdot h(z).
\end{align}
For any function $h$, $\sup_t S_A(h)(t) \leq \sup_t h(t)$, so $\|S_A\|\leq 1$.
Similarly, ${\|S_B\|\leq 1}$.
In order to show that $\| L_B\circ L_A\|<1$, we factorize $L_B \circ L_A$ into ${D_B \circ S_B \circ D_A \circ S_A}$. It then suffices to bound $\|D_A\|$, $\|D_B\|$ or $\|D_B \circ S_B \circ D_A\|$ away from $1$, since by the definition of the operator norm and using that $\|S_A\|,\|S_B\|\leq 1$, we have
\[
\|L_B\circ L_A\| \leq \|D_B \circ S_B \circ D_A\| \leq \|D_A\|\cdot \|D_B\|.
\]
The proof of the lemma will be divided into two cases, depending on whether $I_A(\lambda/2)=I_B(\lambda/2)=1$ or not.
\begin{description}
\item[Case 1: $I_A(\lambda/2)<1 \;\mathrm{ or }\; I_B(\lambda/2)<1$]
\hfill
\\
Assume without loss of generality that $I_A(\lambda/2)<1$. Then $I_A(z)<1$ for all $z$. By \cref{lemma:conditionalshift}, $I_A$ is a continuous function on a closed interval, so it attains its supremum $\tau$, which must be less than $1$. Thus $ \|D_A\| \leq \tau <1$.

\item[Case 2: $I_A(\lambda/2)=I_B(\lambda/2)=1$]
\hfill
\\
$ {D_B \circ S_B \circ D_A}$ is an integral operator with kernel given by $I_B(s)I_A(t)\kappa^s_B(t)$. In order to show that this integral operator has norm less than $1$, we will bound the integral of its kernel along the line $s = z$, where $z \in (-\lambda/2,\lambda/2]$ is arbitrary but fixed. (We do not need to consider the case $z = -\lambda/2$, since $I_B(-\lambda/2)<1$.)

We have a good upper bound on $I_A$ and $I_B$ on any closed set not containing $z=\lambda/2$. In particular, on $[-\lambda/2,0]$. We therefore bound the total mass of $\kappa_B^z$ on $[0,\lambda/2]$:
Since $I_B(\lambda/2) = 1$, we know that $\lint  \rho_B^z(t)dt \to \infty$ as $z\to\lambda/2$.
 But $\int_0^{\lambda/2} \rho_B^z(t) dt \leq q\lambda^{q-1}$ for any $z$, so $\int_0^{\lambda/2} \kappa_B^z(t) dt$ must vanish as $z\to \lambda/2$. Hence there exists $0<\theta<\lambda/2$ such that  for all  $z>\theta$,
\begin{equation}
\int\limits_{-z}^{0}\kappa_B^z(t) dt > 1/2.
\label{bound:massnearsingularity}
\end{equation}
By (ii) of \cref{lemma:conditionalshift}, we can find $\delta>0$, such that when $t\leq \theta$ we have
\begin{equation}
I_A(t)< 1-\delta\quad\textrm{ and }\quad I_B(t)< 1-\delta.
\label{bound:Ileft}
\end{equation}
Then, for $-\lambda/2 \leq z\leq \theta$, we apply the bound to $I_B$ to get
\[
\lint I_B(z)I_A(t)\kappa_B^z(t)dt  \leq I_B(z) \cdot \lint \kappa_B^z(t)dt \overset{(\ref{bound:Ileft})}{<} 1-\delta,
\]
while for $\theta < z \leq \lambda/2$ we apply it to $I_A$
\[
\lint I_B(z)I_A(t)\kappa_B^z(t)dt
\overset{(\ref{bound:Ileft})}{<} \int\limits_{-z}^{0} (1-\delta) \kappa_B^z(t)dt + \int\limits_{0}^{\lambda/2}\kappa_B^z(t)dt 
\overset{(\ref{bound:massnearsingularity})}{<} 1-\delta/2.
\]
This means that the weight along each line of $D_B \circ S_B \circ D_A$ is at most $1-\delta$ when $z\leq \theta$, and at most $1-\delta/2$ otherwise. So in either case, ${ \|D_B \circ S_B \circ D_A\| \leq 1-\delta/2}$.
\end{description}
We conclude that $\| L_B \circ L_A\|<\max(\tau,1-\delta/2)<1$.
\end{proof}

\subsection{The game finishes in finite time }

\noindent Armed with \cref{locallipschitz,lemma:conditionalshift,normbound}, we can proceed with the proof of \cref{finitegame}.

\vspace{1em}

\noindent \begin{proof}[Proof of \Cref{finitegame}]
\label{finitegameproof}
We will use that the positive operator $L:=L_B \circ L_A$ is a contraction (i.e.\ $\|L\|<1$) to construct suitable $\fcn_t$. Let $\mathbf{1}:\Lambda\to \R $ denote the function with constant value $1$. We define the functions $\fcn_t$ (for any $ t\in[0, 2\lambda]$ and some large constants $K,m>0$ to be determined later) by
\begin{equation}
\fcn_t:= K \exp(mt) \cdot\sum_{k=0}^\infty L^{k}\mathbf{1}.
\label{def:Psi}
\end{equation}
By \cref{normbound}, $\|L\|<1$, so the above series is absolutely convergent, whence for all $z\in \Lambda$
\begin{equation}
1\leq \frac{\fcn_t(z)}{K \exp(mt)} \leq \frac{1}{1-\|L\|}.
\label{fcn-norm}
\end{equation}
Because the series is absolutely convergent, we can apply the operator $L$ to $\fcn_t$ by applying it term-wise to the sum in \cref{def:Psi}.
\begin{align}
L \fcn_t(z) = K  \exp(mt) \cdot  \sum_{k=1}^\infty L^{k}\mathbf{1}(z) = \fcn_t(z)-K  \exp(mt).
\label{def:TTonLambda}
\end{align}
Crucially, $L \fcn_t$ is less than $\fcn_t$, and with a sizeable margin.
We will do induction on even $k$ to establish the main claim.
Fix some ${-\lambda/2\leq z \leq \lambda/2}$ and $0\leq t\leq 2\lambda$, and consider the first two moves of the game conditioned on $f_A(\phi)=z$.

Let $u$ be Alice's optimal move from the root (if such a move exists, the result holds trivially otherwise), and $v_0$ the vertex she expects Bob to move to after that. Assume there are $n$ of Bob's move options from $u$ that are sub-optimal w.r.t. $f_A$ but within $t$ of being $f_A$-optimal. Let $v_i$, $1\leq i\leq n$, be the vertices those moves lead to, $t_i:=\delta(v_i)$, $\ell_i$ the cost of the edge $(u,v_i)$, and let $f_i := f_A(v_i)$. By assumption, $t_i\leq t$ for all $1\leq i \leq n$. Note that $t_i$, $v_i$, $f_i$ and $n$ are random functions of $z$. 

For the base case $k=2$, the children of the root in $\Delta_t^2(\phi)$ are the $f_A$-optimal moves (by definition), and there is almost surely at most one such move ($u$, say). Its expected number of children, conditioned on any value of $f_A(u)$, is at most $1+\lambda^q$, so $R^2_t \leq 2+\lambda^q $.
We let $K=2(2+\lambda^q)$, so that $R^2_t \leq K/2< \fcn_t$, establishing the base case.

Next, assume $R_s^k<\fcn_s$ for some even $k > 2$ and all $0\leq s\leq 2\lambda$. We want to show that $R_t^{k+2}<\fcn_t$ as well, and to do that we will bound the expected size of $\Delta_t^{k+2}$.
The tree $\Delta_t^{k+2}$ can be written as an edge-disjoint union of copies of $\Delta$ in the following way:
\begin{equation}
\Delta_t^{k+2}(\phi) = \Delta_t^2(\phi) \cup  \Delta_t^k(v_0) \cup \bigcup_{i=1}^n \Delta_{t-t_i}^k(v_i).
\label{treeinclusion}
\end{equation}
Note that the trees $\Delta_t^k(v_0),\Delta_{t-t_i}^k(v_i)$ are independent conditional on $f_A(v_0),f_A(v_1),$ $\ldots, f_A(v_n)$.
We already have a bound for $\Delta_t^2(\phi)$, and we continue by bounding the conditional expected sizes of $\Delta_t^k(v_0)$ and $\bigcup_{i=1}^n \Delta_{t-t_i}^k(v_i)$.  For $\Delta_t^k(v_0)$, by the definition of $R_t^k$ and \cref{def:TA},
\begin{align}
\E\big[|\Delta_t^k(v_0)|\big|f_A(\phi)=z\big] &= L (R_t^k)(z) \nonumber
\\
&<  L (\fcn_t)(z) \textrm{, since $L$ is a positive operator}\nonumber
\\
&\overset{(\ref{def:TTonLambda})}{=} \fcn_t(z) - K \exp(mt). \label{bound:optimaltree}
\end{align}
Next, we bound the expected size of the union of the trees $\Delta_{t-t_i}^k(v_i)$. To do this we condition first on the random variables $n$, $f_i$ and $t_i$ and then on the event $f_A(\phi)=z$, so that the first conditional expectation is itself a random variable.
\begin{align}
\E\Big[\big|\bigcup \limits_{i=1}^n \Delta_{t-t_i}^k(v_i)\big|\Big| f_A(\phi)=z \Big]
&=\E \bigg[\E\Big[\big|\bigcup \limits_{i=1}^n \Delta_{t-t_i}^k(v_i)\big| \Big| n, t_i, f_i, 1\leq i\leq n \Big] \bigg|f_A(\phi)=z \bigg] \nonumber
\\
&=\E\Big[\sum_{i=1}^n R_{t-t_i}^k(f_i)\Big| f_A(\phi)=z\Big], \label{sumoftrees}
\end{align}
since the subtree rooted in $v_i$, conditioned on $f_A(v_i)$, is independent of $f_A(\phi)$. By the induction hypothesis with $s=t-t_i$,
\begin{equation*}
\sum_{i=1}^n R_{t-t_i}^k(f_i) \leq \sum_{i=1}^n \fcn_{t-t_i}(f_i).
\end{equation*}
Let $\sigma_A$ be the Poisson random measure generated by $\mu_A$ (i.e.\ the counting measure of the points of the Poisson point process with intensity $\mu_A$). It is a sum of Dirac measures, each corresponding to a point in the $\ell f$-square. Among these points, $(\ell_i,f_i)$, $1\leq i \leq n$, are exactly those that lie in the diagonal strip \[D:=\{(\ell,f):z<\ell-f\leq z+t\}\]
Note that for any bounded $\mu_A$-measurable function $h$, we have that $\E[\int \!  h\,  d \sigma_A]=\int \! h \, d\mu_A$ (which can be seen by approximating $h$ by simple functions).
The expression (\ref{sumoftrees}) is then at most
\begin{align}
\E\Big[\sum_{i=1}^n \fcn_{t-t_i}(f_i)\Big|f_A(\phi)=z\Big] &= \E\Big[\iint_D\fcn_{z+t-l+f}(f) \,d\sigma_A(\ell,f)\Big]  \nonumber 
\\
=\iint_D \fcn_{z+t-\ell+f}(f) \,d\mu_A(\ell,f) 
&\overset{\ref{fcn-norm}}{\leq} \iint_D\frac{ K \exp(m(z+t-\ell+f))}{1-\|L\|}\,d\mu_A(\ell,f)  \label{intdiagstrip}
\end{align}
The integrand is constant along diagonals $\ell-f=x$ for fixed $x\in(z,z+t]$. Recall that the one-dimensional measure of such  a diagonal is $J_A^x$. Integrating along these diagonals first, we see that
\begin{align}
(\ref{intdiagstrip}) &\leq \frac{K}{1-\|L\|} \cdot \int\limits_{z}^{z+t} J_A^x   \exp(m(z+t-x)) dx \nonumber
\\
&\overset{(\ref{bound:Jz})}{\leq}  \frac{K\exp(mt)}{1-\|L\|} \cdot  \int\limits_{z}^{z+t} \alpha\lambda^q\cdot\big[(x+\lambda/2)^{q-1}\!+|x-\lambda/2|^{q-1}\big]\exp(m(z-x)) dx \nonumber
\\
&\leq K\exp(mt)\cdot \varepsilon_m,
\label{bound:suboptimaltrees}
\end{align}
for some $\varepsilon_m$ which goes to zero as $m \to \infty$, and does not depend on $k,t$ or $z$. We now have a bound on the expected size of each term in the right hand side of \cref{treeinclusion}. The bounds from \cref{bound:suboptimaltrees,bound:optimaltree} give that
\begin{align}
R^{k+2}_t(z)
&=\E\big[|\Delta_t^{2}(\phi)|
+|\Delta_t^{k}(v)|
+\sum_i |\Delta_{t-t_i}^{k}(v_i)|\big|f_A(\phi)=z\big] \nonumber
\\
&<\Big(K/2\Big)
+\Big(\fcn_t(z)-K\exp(mt)\Big)
+\Big(K\exp(mt) \varepsilon_m\Big). \label{bound:largetree}
\end{align}
Pick $m$ large enough that $\varepsilon_m <1/2$. The expression (\ref{bound:largetree}) is then at most $\fcn_t(z)$, completing the inductive step. Hence $R_t^k \leq \fcn_t$ for all even $k$ and all ${t\in [0, 2\lambda]}$.
\end{proof}

\vspace{1em}

\noindent\begin{proof}[Proof of \cref{prop:unique valuation}]
By \cref{reasonable-game-path}, the game path $P$ is $(\phi,2\lambda)$-reasonable, and is therefore contained in the tree $\Delta_{2\lambda}(\phi)$ of \emph{all} $(\phi,2\lambda)$-reasonable paths. By \cref{finitegame}, $\Delta_{2\lambda}(\phi)$ is almost surely finite, and hence the game finishes after finitely many steps. Thus
$P$ is the finite path  $\phi = u_0 \to u_1 \to \ldots \to u_N$ for some $u_i$'s. For $1\leq i \leq N$, let $\ell_i = \ell(v_{i-1},v_i)$.
Let $S$ be the total payoff for Alice. (The total payoff for Bob is then $-S$.)
Alice pays $\ell_1+\ell_3+\ldots$ to Bob, and Bob pays $\ell_2+\ell_4+\ldots$ to Alice, until one player decides to pay $\lambda/2$ and quit the game. Thus
\[S = -\ell_1+\ell_2-\ell_3 \ldots \pm \ell_N \mp \lambda/2 \]
where the $\pm$-sign depends on whether Alice or Bob is the one to quit, i.e.\ whether $N$ is even or odd. 

\vspace{-1em}
\begin{claim}
\(f_A(\phi) \geq -S\)

\end{claim}
\noindent \begin{proof}[Proof of claim] 
Recall that  $f_A(u_{i-1}) \leq \ell_i-f_A(u_i)$, with equality if $u_{i-1}\to u_i$ is $f_A$-optimal (which is always the case if $i$ is odd).
Using these inequalities along the game path $P$ gives us
\begin{align*}
f_A(\phi)
= \ell_1-f_A(u_1)
\geq \ell_1-\ell_2+f_A(u_2)
= \ldots
\geq \sum_{i=1}^N (-1)^{i+1} \ell_i +(-1)^N  f_A(u_N) .
\end{align*}

If $N$ is even, then Alice is the one that quits, which she only would have done if $f_A(u_N)=\lambda/2$. 
In that case, \[S =  -\lambda/2 + \sum_{i=1}^N (-1)^{i} \ell_i  = -f_A(u_N)+ \sum_{i=1}^N (-1)^{i} \ell_i .\]
If on the other hand $N$ is odd, then Bob is the one that quits, and (like for every vertex)  $f_A(u_N) \leq \lambda/2$. Hence \[S =  \lambda/2 + \sum_{i=1}^N (-1)^{i} \ell_i  \geq f_A(u_N)+ \sum_{i=1}^N (-1)^{i} \ell_i .\]
In either case $-S \leq \sum_{i=1}^N (-1)^{i+1} \ell_i +(-1)^N  f_A(u_N)  \leq f_A(\phi)$.
\end{proof}

By symmetry (reversing the roles of Alice and Bob in the proof, and $\delta$ instead measuring how far Alice deviates from what is $f_B$-optimal) we also have that  \(-f_B(\phi) \geq S\).
Thus $f_A(\phi) \geq -S \geq f_B(\phi)$. But by the choice of $f_A$ and $f_B$, we know that $f_A(\phi) \leq f_B(\phi)$, so we have that $f_A(\phi) = f_B(\phi)$.
For any other $u \in V(T^q_\lambda)$, the subtree rooted in $u$ has the same distribution as the whole $T^q_\lambda$, so a similar argument gives that $f_A(u) = f_B(u)$. Hence $f_A = f_B$. Since $f_A$ and $f_B$ are the maximum and minimum, respectively, in the lattice ordering of all valuations, this implies that the valuation is unique.
\end{proof}

\vspace{1em}
\noindent We end the paper by giving the deferred proof of  \cref{locallipschitz}.

\noindent\begin{proof}[Proof of \cref{locallipschitz}]
\label{proof:locallipschitz}
Let $\mathcal{G}$ be the class consisting of all non-increasing functions $G:\Lambda\to(0,1]$ which satisfy $G(-\lambda/2)=1$. 
Recall that $V$ is the non-linear operator defined by
\[
V(G)(z):=\exp\Big(-\lint q(z+t)_+^{q-1}G(t)dt \Big).
\]
We will find the derivative of $V(G)$ for any $G\in \mathcal{G}$, and then show that $F_A,F_B\in \mathcal{G}$. Since $F_A=V(F_B)$ and $F_B=V(F_A)$\cite[p.1077]{repsym}, this will give us the derivatives $F_A'$ and $F_B'$.
\vspace{-1em}
\begin{claim}
$V(\mathcal{G})\subseteq \mathcal{G}$.
\end{claim}
\noindent\begin{proof}
For any function $G:\Lambda\to\R$, we have that $V(G)(-\lambda/2) =1$, because
\[V(G)(-\lambda/2) := \exp\left(-\lint q(-\lambda/2+t)_+^{q-1} G(t)dt\right)\]
and $(-\lambda/2+t)_+^{q-1}$ vanishes for all $t\in\Lambda$.

Now, pick a $G\in \mathcal{G}$. Note first that $V(G)$ is non-increasing since $G$ is non-increasing and positive. Furthermore, $G\geq 0$ implies $V(G)(z) \leq 1$, and similarly $G\leq 1$ implies
\[V(G)(z)\geq \exp\left(-\lint q(z+t)_+^{q-1}dt)\right) \geq \exp(-\lambda^q)> 0.\]
Thus $V(G)\in \mathcal{G}$, and the claim follows.
\end{proof}
\vspace{-1em}
\begin{claim}
$F_A,F_B\in \mathcal{G}$
\end{claim}
\noindent \begin{proof}
First, note that $F_A$ and $F_B$ are non-increasing by definition. Next, recall that $F_A=V(F_B)$, $F_B=V(F_A)$ and that $V(G)(-\lambda/2)=1$ for any real-valued function on $\Lambda$, whence $F_A(-\lambda/2)=F_B(-\lambda/2)=1$.
\end{proof}

Since any $G\in \mathcal{G}$ is bounded and monotone, it has bounded variation. We can therefore integrate with respect to the measure $dG$, in the sense of a Riemann-Stieltjes integral. However, Riemann-Stieltjes integration is usually defined for non-decreasing functions rather than the non-increasing function $G$ here, and we must be careful with how Riemann-Stieltjes treats the end points of $\Lambda$.
We therefore let $\tilde G$ be defined by $\tilde G:= 1-G$ on $[-\lambda/2,\lambda/2)$ and $\tilde G(\lambda/2):=1$, and work with $d\tilde G$ rather than $dG$.
\vspace{-1em}
\begin{claim}
For any $G\in \mathcal{G}$, $V(G$) is differentiable on the interior of $\Lambda$, with derivative given by
\begin{equation}
\frac{d}{dz}V(G)(z) = V(G)(z) \cdot \lint q(z+t)_+^{q-1}  d\tilde G(t).
\label{guess:F'}
\end{equation}
\end{claim}
\noindent\begin{proof}
To verify \cref{guess:F'}, start by integrating $\lint  q(z+t)_+^{q-1}  d\tilde G(t)$ from $z=-\lambda/2$ to $x$ (for some $x$ with $|x|< \lambda/2$):
\begin{align*}
\int\limits_{-\lambda/2}^{x}\int\limits_{-\lambda/2}^{\lambda/2}  q(z+t)_+^{q-1}  d\tilde G(t)dz
&= \iint\limits_{\substack{-\lambda/2\leq t \leq \lambda/2, \\ -\lambda/2\leq s-t\leq x}} qs_+^{q-1}d \tilde G(t)ds
\\
=\int_0^{\lambda/2-x} \!\!\!qs^{q-1}G(x+s) ds
&=- \ln\big(V(G)(x)\big)
\end{align*}
By the fundamental theorem of calculus, $\ln(V(G)(z))$ is differentiable, with derivative given by $\frac{d}{dz} \ln(V(G)(z)) =- \lint   q(z+t)_+^{q-1}  d \tilde G(t)$. This implies that $V(G)$ is also differentiable, with derivative given by \cref{guess:F'}, proving the claim.
\end{proof}
\vspace{-1em}
\begin{claim}
Let the function $g$ on $\Lambda$ be defined by 
\begin{equation}
g(z) := (\lambda/2-|z|)^{q-1}.
\label{def:g}
\end{equation}
Then there exists a constant $a>0$ such that if $G\in V(\mathcal{G})$ satisfies $-G' \leq a g$, then $-(V(G))'\leq a g$.
\end{claim}
\noindent\begin{proof}
We need to calculate (and then estimate) $\frac{d}{dz}V(G)(z)$. Since $G\in V(\mathcal{G})$, $G$ is differentiable, and therefore $d\tilde G(t) = -G'(t)dt$ for $t$ in the interior of $\Lambda$. However, $\tilde G$ also has a point mass at $\lambda/2$, so for any $t\in \Lambda$ we have that \[d\tilde G(t) = -G'(t)dt+G(\lambda/2)d\delta_{\lambda/2}(t),\]
where $\delta_x$ is a Dirac measure at $x$. Substituting this expression for $d\tilde G$ in \cref{guess:F'},
\begin{align}
\frac{d}{dz}V(G)(z) 
&=-V(G)(z) \cdot \bigg(q(\lambda/2+z)^{q-1} G(\lambda/2) - \lint   q(z+t)_+^{q-1}  G'(t) dt \bigg).
\label{def:F'}
\end{align}
The integrand on the right hand side of \cref{def:F'} has a singularity at $t=-z$, while the function $g$ has a singularity at $t=\lambda/2$. For some positive parameter $r<\min(\lambda/4,2^{-4/q})$, we will deal separately with two cases: when these singularities are within $2r$ of each other, and when they are further apart.
We will establish that the following inequality holds in both cases:
\begin{equation}
-\frac{d}{dz}V(G)(z) < g(z) \cdot \Big( 2q +4 a r^q\Big)+qr^{q-1},\textrm{ for all } z
\label{bound:case1case2}
\end{equation}
from which it follows that $-\frac{d}{dz}V(G)(z)<a\cdot g(z)$ for all $z$ by picking ${a >\max(8q,\lambda/r)}$. 
\begin{description}
\item[Case 1: $-\lambda/2 \leq z \leq -\lambda/2+2r$] \hfill \\
We apply the bound $-G'(t)\leq a \cdot g(t)$, and use that the resulting integrand is symmetric around $t=-z/2+\lambda/4$:
\begin{align}
-\lint  G'(t)\cdot q(z+t)_+^{q-1} dt
&\leq a\lint   q(\lambda/2-t)^{q-1}(z+t)_+^{q-1} dt \nonumber
\\
&\leq 2a (z/2+\lambda/4)^{q-1} \cdot \int_{0}^{z/2+\lambda/4} q s^{q-1}ds \nonumber
\\
&= 2^{2-2q} a (z+\lambda/2)^{2q-1} \label{bound:2q-1}
\\
&< 4ar^q g(z). \label{bound:case2-integral}
\end{align}
We will later be using the tighter bound in \cref{bound:2q-1}, but for now \cref{bound:case2-integral} suffices. Again using \cref{def:F'}, this gives a bound on $\frac{d}{dz}V(G)(z)$:
\[
-\frac{d}{dz}V(G)(z) \leq G(z) \cdot \bigg(G(\lambda/2) \cdot q(\lambda/2+z)_+^{q-1}+ 4ar^q g(z)\bigg) 
\leq g(z) \cdot \Big( q + 4ar^q \Big)
\]
which is less than the bound from \cref{bound:case1case2}.

\item[Case 2: $-\lambda/2+2r \leq z \leq \lambda/2$] \hfill \\
We use the bound $-G'(t)\leq a \cdot g(t)$ for $-z<t<-z+r$.
\begin{equation}
-\! \lint  G'(t)\cdot q(z+t)_+^{q-1} dt \leq
 \!\!\int\limits_{-z}^{-z+r}  \!\!ag(t)\cdot q(z+t)^{q-1} dt \; 
-  \!\!\int\limits_{-z+r}^{\lambda/2}   \!\!G'(t)\cdot q(z+t)^{q-1} dt
\label{bound:case1}
\end{equation}
If $g(t)$ is larger than $g(-z)$, for $-z\leq t \leq -z+r$, it can be at most twice as large, since $g$ is increasing fastest at $\lambda/2-2r$ and $g(\lambda/2-r) \leq 2 g(\lambda/2-2r)$.
Hence the first integral on the right hand side of \cref{bound:case1} is at most
\begin{equation}
\int_{-z}^{-z+r}  2g(-z) \cdot q (z+t)^{q-1}dt \leq 2g(z)\cdot r^q, \label{bound:case1-1st-integral}
\end{equation}
while second integral on the right hand side of \cref{bound:case1} is at most
\begin{equation}
-\int_{-z+r}^{\lambda/2}   G'(t)\cdot  qr^{q-1} dt \leq qr^{q-1},
\label{bound:case1-2nd-integral}
\end{equation}
since $q(z+t)^{q-1}$ is a decreasing function in $t$. Putting \cref{bound:case1-1st-integral,bound:case1-2nd-integral} together with \cref{def:F'}, this gives us that $-\frac{d}{dz}V(G)(z)$ is at most
\[
V(G)(z) \cdot \bigg(G(\lambda/2) q(\lambda/2+z)^{q-1}+2a\cdot g(z) r^q +  qr^{q-1}\bigg)
\leq 2 g(z) \cdot \Big( q +ar^q\Big)+qr^{q-1},
\]
which is also less than the bound from \cref{bound:case1case2}.
\end{description}
\vspace{-2.03em}
\end{proof}
\begin{claim}
$-F_A',-F_B'\leq ag$  (ineq. (\ref{bound:F'}) in the statement of the lemma)
\end{claim}
\noindent \begin{proof}
Note first that $F_A,F_B\in V(\mathcal{G})$, whence they are differentiable by a previous claim.
Let $G_1(z):=1$ for all $z\in \Lambda$, and $G_{k+1}:=V(G_{k})$. Then $G_1\in \mathcal{G}$, and by induction $G_k\in \mathcal{G}$ for all $k\geq 1$. We know by \cite[p.1077]{repsym} that for any $z\in \Lambda$, $G_{2k}(z)\nearrow F_A(z)$, and similarly $G_{2k+1}(z)\searrow F_B(z)$.

$\Lambda$ is compact and $F_A$ is continuous, so by Dini's theorem $G_{2k}\to F_A$ uniformly.
But since $F_A$ is differentiable, uniform convergence of $G_{2k}$ implies $G'_{2k}\to F_A'$. Similarly, $G'_{2k+1}\to F_B'$ as $k\to \infty$.

Noting that $-G_1'(z)=0< ag$, and $-G'_k(z)<ag\Rightarrow -G'_{k+1}(z)<ag$, by induction $-G'_k(z)<ag$ for all $k$. Since $G_{2k}'\to F_A'$, we have that $-F_A'\leq ag$, and similarly $-F_B'\leq ag$.
\end{proof} 
 
The next step is to show that $\lint  \rho^z_A(t)dt$ is continuous in $z$ on $(-\lambda/2,\lambda/2)$. (Recall that $\rho^z_A(t):= -q(t+z)_+^{q-1}F_A'(t)$.) It suffices to show that it is continuous on any closed subinterval ${I \subset (-\lambda/2,\lambda/2)}$. Since $|F'_A|$ is bounded by $K_I:= \sup_{t\in I} ag(t)< \infty$ on $I$, for any $x,y \in I$ we have
\begin{equation}
|F_A(x)-F_A(y)|<K_I\cdot |x-y|
\label{bound:lipschitz}
\end{equation}
We will let $\varepsilon := \sqrt{|x-y|}\to 0$. Suppose (without loss of generality) that $x<y$ and ${[x-\eps,y+\eps]\subseteq I}$.
We estimate the difference
\begin{align*}
&\quad \quad \left|\int\limits_{-x}^{\lambda/2} \!q(x+t)^{q-1}F'_A(t)dt-\int\limits_{-y}^{\lambda/2} \!q(y+t)^{q-1}F'_A(t)dt\right| 
\\
&\overset{(\ref{bound:lipschitz})}{\leq} 2\left| \int\limits_{-x}^{-x+\varepsilon} \!K_I \cdot q(x+t)^{q-1} dt \right| 
+q\left|\int\limits_{-y+\varepsilon}^{\lambda/2}  \Big((x+t)^{q-1}-(y+t)^{q-1}\Big) F'_A(t)dt \right|
\\
&\leq 2K_I \eps^q+ 2\underbrace{|x-y|}_{\leq \eps^2} \cdot q(1-q)\cdot \left|\int\limits_{-y+\varepsilon}^{\lambda/2}  \underbrace{(y+t)^{q-2} }_{\leq \eps^{q-2}}F'_A(t)dt \right| = O(\varepsilon^q).
\end{align*}
Hence $\lint  \rho_A^z(t)dt$ is continuous in $z$, and so is $F'_B$.

To establish the bound (\ref{bound:diag}) for $\lint  \rho_A^z(t)dt$, we use that $-F'_A \leq a g$. Then $F_A$ satisfies the conditions necessary for \cref{bound:2q-1} to hold for $z$ near $-\lambda/2$ with $G=F_A$. For other $z$, note that the integrand is at most $-F'_B(z)$, for which the weaker bound $a g$ suffices. In other words, for some constant $b$ and any $-\lambda/2<z<\lambda/2$, we have that
\[
\lint  \rho^z_A(t)dt \leq b \max\big((\lambda/2-z)^{q-1},(z+\lambda/2)^{2q-1}\big).
\]
Finally, for $z>\lambda/2$, note that ${(z+t)^{q-1}\leq(z-\lambda/2)^{q-1}}$, whence ${\lint \rho^z_A(t)dt }$ is at most $ {q(z-\lambda/2)^{q-1}}$. Setting $\alpha = \max(a,b,q)$ gives the desired result.
\end{proof}

\bibliographystyle{plain} 
 
\newpage
\end{document}